\documentclass[11pt]{article}

\usepackage[leqno]{amsmath}
\usepackage{amssymb} 

\usepackage{graphicx} 
\usepackage[applemac]{inputenc}

\newcommand{\sca}[2]{\langle #1, #2\rangle}

\newcommand{\R}{\mathbb R}
\renewcommand{\P}{\mathbb P}
\newcommand{\C}{\mathbb C}

\let\p=\partial

\let\lb=\lambda 
 
\let\de=\delta

\let\vf=\varphi

\let\sm=\setminus
\let\ss=\subset

\let\w=\omega

\let\i=\infty
\let\t=\tau
\let\hat=\widehat 
\let\bar=\overline

\let\D=\Delta
\let\tr=\textrm

\begin{document}

\centerline{\bf  \large A hypersurface containing the support of a Radon transform} 
\smallskip
\centerline{\bf \large must be an ellipsoid. I}

\medskip
\centerline{ Jan Boman}
  
\medskip
\centerline{(2019-07-30)}

\begin{abstract}
If the Radon transform of a compactly supported distribution $f \ne 0$ in $\R^n$ is supported on the set of tangent planes to the boundary $\p D$ of a bounded convex domain $D$, then $\p D$ must be an ellipsoid. As a corollary of this result we get a new proof of a recent theorem of  Koldobsky,   Merkurjev,  and Yaskin, which settled a special case of a conjecture of Arnold that was motivated by a famous lemma of Newton. 
\end{abstract}
 
\bigskip
\noindent
{\bf 1. Introduction.} 
Define a function $f_0$ in the plane by 
\begin{equation*}
f_0(x) =  \frac 1{\pi} \, \frac 1{\sqrt{1 - |x|^2}},  \quad \tr{for $x = (x_1, x_2) \in \R^2$, $|x| < 1$} , 
\end{equation*}
and $f_0(x) = 0$ for $|x| > 1$. 
An easy calculation shows that the Radon transform of $f_0$ satisfies 
\begin{equation*}
R f_0(\w, p) = \int_{x \cdot \w = p} f_0 \, ds = 1 \quad \tr{for $|p| < 1$, \  $\w \in S^1$}  , 
\end{equation*}
and obviously $R f_0(\w, p) = 0$ for $|p| \ge 1$. 
Define the distribution $f$ by 
\begin{equation*}
f = \D f_0 = \p_{x_1}^2 f_0 + \p_{x_2}^2 f_0 . 
\end{equation*}
The  well known formula $R(\D h)(\w, p) = \p_p^2 Rh(\w, p)$ with $h = f_0$ now
shows  that  
\begin{equation*}
R f(\w, p) = 0  \quad \tr{for $|p| < 1$ and all $\w$} . 
\end{equation*} 
In other words,  $R f$ is a distribution on the manifold of lines in the plane that vanishes in the open set of lines that intersect the open unit disk. Since $R f$ obviously vanishes on the open set of lines that are disjoint from the closed disk, it follows that the distribution $R f$ is 
supported on the set of lines that are tangent to the circle. 

By means of an affine transformation it is easy to construct a similar example where the circle is replaced by an arbitrary ellipse.  For an arbitrary ellipsoidal domain $D \ss \R^n$, $n>2$,  it  is also easy to construct examples of distributions $f$ supported  in $\bar D$ such that the Radon transform $R f$ is supported on the set of tangent planes to the boundary of $D$. 
However, surprisingly, for other convex domains than ellipsoids such distributions do not exist. 

Since we will consider arbitrary convex  – not necessarily smooth – domains, we have to replace the notion of tangent plane by supporting plane. A supporting plane for $D$ is a hyperplane $L$ such that $L \cap \bar D$ is non-empty and one of the components of $\R^n \sm L$ is disjoint from $D$. 

\smallskip 
\noindent
{\bf Theorem 1.}
{\sl Let $D$ be an open, convex, bounded and symmetric (that is $D = - D$) subset of $\R^n$, $n\ge 2$, 
with  boundary $\p D$. If there exists a distribution $f \ne 0 $ with support in $\bar D$ such that the Radon transform of $f$ is supported on the set of supporting planes for $D$, then $\p D$ must be an ellipsoid. 
}

The more general case when $D$ is not assumed to be symmetric turned out to require  different arguments from those given here. This case will therefore be treated in a forthcoming article. 

\smallskip 
\noindent
{\it Remark.}
Our arguments prove in fact a stronger statement, Theorem 3, which is local in $\w$ but global in $p$; see Section 5. 

Denote by $V(\w, p)$ the volume of one of the components of $D\sm L(\w, p)$, where $D\ss \R^n$ is a convex, bounded domain and $L(\w, p)$ is the hyperplane $x \cdot \w = p$ that we assume intersects $D$.  A famous conjecture of Arnold (Problem 1987-14 in {\it Arnold's Problems}, \cite{A}) asserts that if $V(\w, p)$ is an algebraic function, then $n$ must be odd and $D$ must be an ellipsoid. The case $n$ even has been settled long ago by Vassiliev \cite{V1}, see also \cite{V2}. For odd $n$ the question is still open. 
The special case when $n$ is odd and $p\mapsto V(\w, p)$ is assumed to be a polynomial function of degree $\le N$ for all $\w$ and some $N$ has also been studied and was   settled recently by Koldobsky,   Merkurjev,  and Yaskin  for domains $D$ with smooth boundary, \cite{KMY}; see also \cite{AG1}.   If the domain $D$ is symmetric, then this question is answered by Theorem~1.  In \cite{AG2} the case when $p\mapsto V(\w, p)$ is algebraic and satisfies a certain additional condition is reduced to the case $p\mapsto V(\w, p)$ is polynomial. 

\smallskip
{\noindent}
{\bf Corollary 1.} 
{\sl Let $D \ss \R^n$, $n \ge 2$,  be as in Theorem 1, and assume that there exists a number $N$ such that 
$p\mapsto V(\w, p)$ is a polynomial of degree $\le N$ for all $\w \in S^{n-1}$. 
Then $\p D$ must be an ellipsoid. 
}

{\noindent}
{\it Proof.}
Let $\chi_D$ be the characteristic function of $D$ and choose an integer $k$ such that $2k > N$.  The assumption implies that 
$\p_p^{2k} (R \chi_D)(\w, p) = 0$ for all $p$ in the interval 
\begin{equation*}
\inf\{x \cdot \w;\, x \in D\} < p < \sup\{x \cdot \w;\, x \in D\}   , 
\end{equation*}
and obviously $(R \chi_D)(\w, p) = 0$ for all other $p$. This shows that the distribution $\p_p^{2k} R \chi_D$ must be supported on the set of supporting planes to $\p D$, a hypersurface in the manifold of hyperplanes in $\R^n$. 
Define the distribution $f$ in $\R^n$ by $f = \D^k \chi_D$ where $\D$ denotes the Laplace operator. 
The formula $R(\D^k h)(\w, p) = \p_p^{2k} Rh(\w, p)$ with $h = \chi_D$ now shows that the distribution $R f $  must be supported on the  set of supporting hyperplanes.  By Theorem~1 this implies that $\p D$ is an ellipsoid. 

A somewhat related problem is treated in a recent article by Ilmavirta and Paternain, \cite{IP}. It is proved that the existence of a function in $L^1(D)$, $D \ss \R^n$,  whose X-ray transform (integral over lines) is constant, implies that the boundary of $D$ is a ball. 

In Section 2 we will write down an expression for an arbitrary distribution $g$ on the manifold of hyperplanes that is equal to the Radon transform of some compactly supported distribution and is supported on the submanifold of supporting planes to $\p D$. In Section 3 we will use the description of the range of the Radon transform  to write down the conditions for $g$ to be the Radon transform of a compactly supported distribution $f$. Those conditions will be an infinite number of polynomial identities in the supporting function $\rho(\w)$ for $D$ and the densities $q_j(\w)$ that define the distribution $g$. Thereby the problem is transformed to a purely algebraic question. In Section~4 we analyze the polynomial identities and prove (Theorem~2) that they imply that   the supporting function $\rho(\w)$ must be a quadratic polynomial, which together with the fact that $\rho(\w) > 0$ implies that $\p D$ is an ellipsoid. In Section 5 we finish the proof of Theorem~1 and prove the semi-local version Theorem~3.  An outline of the proof of Theorem~2 is given in Section~4.

\medskip
\noindent
{\bf 2. Distributions on the manifold of hyperplanes.}
As is well known the manifold $\P^n$ of hyperplanes in $\R^n$ can be identified with the manifold $(S^{n-1}\times \R)/(\pm 1)$, the set of pairs $(\w, p) \in S^{n-1}\times \R$, where $(\w, p)$ is identified with $(-\w, -p)$. Thus a function on $\P^n$ can be represented as an even function $g(\w, p) = g(-\w, -p)$ on $S^{n-1}\times \R$. In this article a distribution on $\P^n$ will be a linear form on 
$C_e^{\i}(S^{n-1}\times \R)$, the set of  smooth even functions on $S^{n-1}\times \R$, and a 
locally integrable even function $h(\w, p)$ on  $S^{n-1}\times \R$ will be identified with the distribution 
\begin{equation*}
C_e^{\i}(S^{n-1}\times \R) \ni \vf \mapsto \int_{\R} \int_{S^{n-1}} h(\w, p) \vf(\w, p) d\w\, dp , 
\end{equation*}
where $d\w$ is area measure on $S^{n-1}$. 
Using the standard definition of $R^*$,
\begin{equation*}
R^*\vf(x)(x) = \int_{S^{n-1}}\vf(\w, x \cdot\w) d\w  , 
\end{equation*}
we can then define the Radon transform of the compactly supported distribution $f$ on $\R^n$ as the linear form 
\begin{equation*}
C_e^{\i}(S^{n-1}\times \R)\ni \vf \mapsto \sca f{R^* \vf} . 
\end{equation*}

Let $D$ be a bounded, convex subset of $\R^n$ with  boundary $\p D$. Here we will also assume that $D$ is symmetric with respect to some point, which we may assume to be the origin, so $D = -D$. We shall denote the supporting function of $D$ by $\rho(\w)$, that is 
\begin{equation*}
\rho(\w) = \sup\{x\cdot \w;\, x \in D \} . 
\end{equation*}
Since $D$ is symmetric, $\rho$ will be an even function, because 
\begin{align*}
\rho(-\w)  &= \sup\{x\cdot(- \w);\, x \in D \} = \sup\{x\cdot (-\w);\, -x \in D \}   \\ 
& = \sup\{(-x)\cdot (-\w);\, x \in D \} = \rho(\w) . 
\end{align*}
Clearly a hyperplane $x\cdot\w = p$  intersects $D$ if and only if $|p| < \rho(\w)$, and it 
 is a supporting plane to $\p D$ if and only if  $p = \pm \rho(\w)$. 
We shall consider the hypersurface in $\P^n$ that consists of all the supporting planes to $\p D$. Since the origin in $\R^n$ is contained in (the interior of) $D$, none of the supporting planes can contain the origin, hence $\rho(\w) > 0$ for all $\w$.

 A distribution of order $0$  that is supported on the set of supporting planes to $D$ can therefore be represented as 
 \begin{equation*}
 g(\w, p) = q_+(\w) \de(p - \rho(\w)) + q_-(\w)\de(p +\rho(\w))
 \end{equation*} 
 for some measures $q_+(\w)$ and $q_-(\w)$ on $S^{n-1}$; here $\de(\cdot)$ is the Dirac measure at the origin in $\R$.   Since $\rho(-\w) = \rho(\w)$ and $\de(\cdot)$ is even we have 
 \begin{align*}
 g(-\w, -p) & = q_+(-\w) \de(-p - \rho(\w)) + q_-(-\w)\de(-p +\rho(\w))  \\ 
     & = q_+(-\w) \de(p + \rho(\w)) + q_-(-\w)\de(p -\rho(\w))  .  
 \end{align*}
Since $g$ must be even, 
 $g(\w, p) = g(-\w, -p)$, this shows that we must have $q_-(-\w) = q_+(\w)$. 
Denoting $q_+(\w)$ by $q_0(\w)$ we can therefore write 
  \begin{equation}    \label{g1}
 g(\w, p) = q_0(\w) \de(p - \rho(\w)) + q_0(-\w)\de(p +\rho(\w))
 \end{equation} 
 for some measure $q_0(\w)$. 
 
 We next show that we may assume that the distribution $f$ is even, $ f(x) = f(-x)$, which implies that $g = Rf$ is even in $\w$ and $p$ separately.

 \noindent
 {\bf Lemma 1.}
 {\sl 
Assume that there exists a compactly supported distribution $f \ne 0$ such that $R f$ is supported on $p = \pm \rho(\w)$. 
Then there exists an even distribution with the same property. 
}

\noindent
 {\it  Proof.}
 Let $f \ne 0$ be such that  $R f$ is supported on $p = \pm \rho(\w)$. We have to construct an even distribution with the same property. 
It is   clear that the distribution $f(-x)$ has the same property. Hence the even part 
$(f(x) + f(-x))/2$ and the odd part $(f(x) - f(-x))/2$ of $f$ both have the same property. If the even part is different from zero there is nothing more to prove, so we may assume that the odd part  is different from zero.  Set $h(x) = (f(x) - f(-x))/2$. Then $h_1 = \p h/\p x_1$ is an even distribution. It remains to prove that $R h_1$ is supported on  $p = \pm \rho(\w)$. But this follows from the formula $R (\p_{x_1} h)(\w, p) = \w_1 \p_p R h(\w, p)$, which is easily seen by application of the formula 
$\hat{R \vf}(\w, \t) = \hat{\vf}(\t\w)$ to both members.

From now on we will therefore assume that the distribution $f$ is even, which implies that $g(\w, p) = Rf(\w, p)$ is even in $\w$ and $p$ separately. This implies that the measure $q_0$ in \eqref{g1} must be even,  so we may write 
 \begin{equation}    \label{g2}
 g(\w, p) = q_0(\w) \big(\de(p - \rho(\w)) +  \de(p +\rho(\w)\big)
 \end{equation} 
 for some $q_0(\w)$.

If the boundary $\p D$ is smooth and hence $\rho(\w)$ is smooth, we can argue 
similarly, using the fact that $\de^{(j)}(\cdot)$ is odd if $j$ is odd and even if $j$ is even, to see that an arbitrary   distribution $g(\w, p)$ that is even in $\w$ and $p$ separately  and is supported on $p = \pm \rho(\w)$ can be written 
\begin{equation}    \label{g}
g(\w, p) = \sum_{j=0}^{m-1}
q_j(\w)  \big( \de^{(j)}(p - \rho(\w)) + (-1)^j  \de^{(j)}(p +\rho(\w)) \big) 
\end{equation}
for some   distributions $q_0(\w), \ldots , q_{m-1}(\w)$ on $S^{n-1}$.  
But if $\rho(\w)$ is not smooth, this is not always true. 
Note that $\de^{(j)}(p \pm \rho(\w))$ should be interpreted as the $j$th distribution derivative of  $\de(p \pm \rho(\w))$ with respect to $p$. 

However, if $g = R f$ for some compactly supported distribution $f$, then we shall see that the representation \eqref{g} is valid and that the distributions $q_j(\w)$ must be continuous functions.

\smallskip
\noindent
{\bf Lemma 2.}
{\sl 
Let $f$ be a compactly supported even distribution in $\R^n$ and let $g = R f$. Assume that $g$ is supported on the set of supporting planes to $D$. Then there exists a number $m$ and continuous functions $q_j(\w)$ such that the distribution $g$ can be written in the form \eqref{g}. 
} 

\noindent
{\it Proof.}
For arbitrary $\w \in S^{n-1}$ define the distribution $R_{\w} f$ on $\R$ by 
\begin{equation*}
\sca{R_{\w} f}{\psi} = \sca f{x \mapsto \psi(x \cdot \w)} \quad \tr{for $\psi \in C^{\i}(\R)$} . 
\end{equation*}
We  note that the map $\w \mapsto R_{\w} f$ must be continuous in the sense that 
$\w \mapsto \sca{R_{\w} f}{\psi}$ is continuous for every test function $\psi \in C^{\i}(\R)$.  
$R f$ can be expressed in terms of  $R_{\w} f$ as follows. If $\vf(\w, p) = \vf_0(\w) \vf_1(p)$, then 
\begin{align}     \label{Rf}
\begin{split}
\sca{R f}{\vf} & = \sca f {R^*\vf} = \sca f{\int_{S^{n-1}} \vf_0(\w) \vf_1(x \cdot \w)} d\w   \\ 
  & =   \int_{S^{n-1}} \vf_0(\w) \sca f{\vf_1(x \cdot \w)} d\w =  \int_{S^{n-1}} \vf_0(\w) \sca{R_{\w} f}{\vf_1} d\w . 
  \end{split}
\end{align}
To prove the second last identity we replace the integrals by Riemann sums and observe that the function $x \mapsto \vf_1(x \cdot \w)$ together with all its derivatives depends continuously on $\w$. 
The formula \eqref{Rf} shows that if  $g = R f$ is supported on the hypersurface 
$p = \pm \rho(\w)$, then  $R_{\w} f$  must be supported on the union of the two points $p = \pm \rho(\w)$ for every $\w$. Hence $R_{\w} f$ can be represented as the right hand side of \eqref{g} for every $\w$. It remains only to prove that all   $q_j(\w)$ are  continuous. It is enough to prove that   $q_j(\w)$ is continuous in some neighborhood of an arbitrary $\w^0 \in S^{n-1}$. If we choose $\psi$ such that $\psi(p) = 0$ in some neighborhood of $-  \rho(\w^0)$ then    
\begin{equation*}
\sca{R_{\w} f}{\psi} = \sum_{j=0}^m  q_j(\w)  \sca{\de^{(j)}(p - \rho(\w))}{\psi(p)} 
   =   \sum_{j=0}^m (-1)^j q_j(\w) \psi^{(j)}(\rho(\w)) . 
\end{equation*}
We have seen that the expression on the right hand side must be a continuous function of $\w$ for every $\psi$. 
Choosing $\psi(p)$ such that $\psi(p) = 1$ in a neighborhood of $ \rho(\w^0)$ shows that $q_0(\w)$ is continuous at $\w^0$. 
Next choosing $\psi(p)$ such that $\psi(p) = p$ in a neighborhood of $ \rho(\w^0)$ shows that $q_1(\w)$ is continuous. Continuing in this way completes the proof. 

\smallskip

Our next goal will be to write down the conditions on $q_j(\w)$ and $\rho(\w)$ for $g(\w, p)$ to belong to the range of the Radon transform.

 \medskip
\noindent
{\bf 3. The range conditions.}
It is well known that a compactly supported $(\w, p)$-even function or distribution $g(\w, p)$ belongs to the range of the Radon transform if and only if the function 
\begin{equation*}
\w = (\w_1, \ldots , \w_n) \mapsto \int_{\R} g(\w, p) p^k dp 
\end{equation*}
is the restriction to the unit sphere of a homogeneous polynomial of degree $k$  in $\w$ for every non-negative integer $k$.

We next compute the moments  $ \int_{\R} g(\w, p) p^k dp$ for the expression \eqref{g}.  
By the definition of $\de^{(j)}$, for any  $a \in \R$,  $ \int_{\R} \de^{(j)}(p - a) p^k dp = 0$ if  $j > k$ and 
\begin{equation*}
 \int_{\R} \de^{(j)}(p - a) p^k dp  
 =  (-1)^j \int_{\R} \de(p - a) \, \p_p^j p^k dp  
  =  (-1)^j \frac{k!}{(k-j)!} a^{k-j} , \quad \tr{if $j \le k$} .
\end{equation*}
Hence if $j \le k$  and $a \ne 0$ we get  
$\int_{\R} \big( \de^{(j)}(p - a)   + (-1)^j   \de^{(j)}(p + a) \big) p^k dp = 0$ if $k$ is odd and 
\begin{equation*}
 \int_{\R} \big( \de^{(j)}(p - a)  +   (-1)^j  \de^{(j)}(p + a)  \big) p^k dp   
 = 2 (-1)^j  \frac{k!}{(k-j)!} a^{k-j}  
\end{equation*} 
if $k$ is even. For arbitrary non-negative integers $k, j$ we define the constant $c_{k, j}$ by $c_{0,0} = 1$ and 
\begin{align}       \label{c1}
\begin{split}
c_{k, j} &= \frac{k!}{(k-j)!} = k (k-1) \ldots (k-j+1) \quad \tr{if  $0 \le j \le k$ and $k \ge 1$},     \\
c_{k, j} &= 0  \quad \tr{if  $j > k$} .   
\end{split}
\end{align}
Note that the second expression   on the first line above  makes sense also for $j > k$ and is equal to zero then, although the first expression does not make sense if $j > k$. For instance, if $j = 2$, then $c_{k,j} = k (k-1)$ for all $k$. 
Then we can now summarize our computations as follows: if $g(\w, p)$ is defined by \eqref{g}, then $\int_{\R} g(\w, p) p^k dp = 0$ if $k$ is odd and 
\begin{equation}    \label{moment}
\int_{\R} g(\w, p) p^k dp = 2 \sum_{j=0}^k c_{k, j} (-1)^j q_j(\w) \rho(\w)^{k-j}, \quad \tr{if $k$ is even} . 
\end{equation}
Thus, for $g(\w, p)$ to be the Radon transform of a compactly supported distribution it is necessary and sufficient that 
\begin{equation}     \label{moments}
\sum_{j=0}^k c_{k, j} (-1)^j q_j(\w) \rho(\w)^{k-j} 
\end{equation}
is equal to the restriction to $S^{n-1}$ of a homogeneous polynomial  for every even $k$. In the next section we will show that those conditions imply that $\rho(\w)^2$ must be a quadratic polynomial.

 The fact that $\rho(\w)^2$ is a quadratic polynomial, combined with the fact that  $\rho(\w)$ is everywhere positive on $S^{n-1}$,  implies that $\p D$, the boundary of the region $D$, is an ellipsoid. Indeed, if  $\rho(\w)$ is also rotationally invariant, $\rho(\w) = c |\w|^2$, then it is obvious that $D$ must  be rotationally invariant, hence a ball. And since $\rho(\w)$ is (strictly) positive we can find an affine transformation $A$ such that $\rho(A\w)^2 = c|\w|^2$ for some $c$.  This implies that $D$ must be an affine image of a ball, hence an ellipsoid.

\smallskip
\noindent
{\bf 4. Analysis of the polynomial identities.}
The purpose of this section is to prove the following purely algebraic result. 
We shall denote the set of restrictions to the unit sphere of   homogeneous polynomials of degree $k$ by $\mathcal P_k$. 
\smallskip
  
\noindent
{\bf Theorem 2.}
{\sl 
Assume that the strictly positive and continuous function $\rho(\w)$ on $S^{n-1}$ and the continuous functions  $q_0, q_1, \ldots , q_{m-1}$, not all zero,  satisfy the infinitely many identities 
\begin{equation}       \label{id1}
 \sum_{j=0}^{m-1} c_{2k,j}  \rho(\w)^{2k-j} q_j(\w)  = p_{2k}(\w) \in \mathcal P_{2k} \quad \tr{for $k = 0, 1, \ldots$ ,}  
\end{equation}
where $c_{k, j}$ is defined by \eqref{c1}. Then $\rho(\w)^2$ is a (not identically vanishing) quadratic polynomial. 
}

In \eqref{id1} we have omitted the factor $(-1)^j$ that occurred in \eqref{moments}, because in the proof of Theorem 1 we can of course apply Theorem 2 to the functions $(-1)^j q_j$. 

For instance, if $m = 3$ the first few equations \eqref{id1} read 
\begin{align}    \label{id1a}
\begin{split}
& q_0  = p_0   \\  
 &  q_0  \rho ^2 + 2 \, q_1  \rho  + 2 \, q_2 = p_2    \\ 
  &   q_0  \rho ^4 + 4 \, q_1  \rho ^3 + 4 \cdot 3 \,  q_2  \rho ^2   = p_4   \\  
  &  q_0  \rho ^6 + 6 \, q_1  \rho ^5 + 6 \cdot 5 \,  q_2  \rho ^4     = p_6 \\  
  &  q_0  \rho ^8 + 8 \, q_1  \rho ^7 + 8 \cdot 7 \,  q_2  \rho ^6     = p_8 . 
  \end{split}
\end{align}

\smallskip

The first step of the proof of  Theorem 2 will be to eliminate the $m$ functions $q_j$ from sets of $m + 1$ of the equations \eqref{id1}. The result is a set of infinitely many polynomial identities in $\rho(\w)^2$ with the polynomials $p_k$ as coefficients, as will be explained in Lemma~4.  
Considering a set of $m$ of those identities as a linear system of equations in the $m$ quantities $\rho^2, \rho^4,  \ldots , \rho^{2m}$  we can solve 
$\rho^2$ as a rational function in the coefficients $p_{2k}$ and hence as a rational function of $\w$, provided the determinant of the corresponding coefficient matrix \eqref{A0} is different from zero. Under the same assumption we can prove rather easily that $\rho^2$ must be a polynomial by considering sets of $m$ such linear systems together. This entails considering the translation operator 
\begin{equation*}
(p_{2k}, p_{2k+2}, \ldots , p_{2k+2m-2}) \mapsto (p_{2k+2}, p_{2k+4}, \ldots , p_{2k+2m}) 
\end{equation*} 
on $m$-vectors of polynomials from 
the infinite sequence  $p_0, p_2, \ldots$.  This operator is given by the matrix $S$ introduced below \eqref{S}. The matrix $S$ has $m$ identical eigenvalues $\rho^2$, hence its determinant is $\rho^{2m}$.  The crucial fact that the matrix \eqref{A0}  is non-singular (Proposition 1) is an easy consequence of the fact that that Jordan canonical form of $S$ consists of just one Jordan block (Lemma 6).

\smallskip
\noindent
{\bf Lemma 3.}
{\sl 
The rank of every $m \times m$ minor of the infinite matrix $c_{2k,j}$, $k =  0, 1, 2 \ldots$, $0 \le j \le m-1$, is maximal, that is, equal to $m$.    
}

\smallskip
\noindent
{\it Proof.}
Introduce the polynomials 
\begin{equation*}
h_0(t) = 1, \quad h_1(t) = t, \quad  h_2(t) = t(t - 1) , \quad \ldots ,  \quad h_j(t) = t (t-1) \ldots (t-j+1) . 
\end{equation*}
Then $c_{2k,j} = h_j(2k)$. Since for any $j$
\begin{equation*}
t^j = h_j(t) + \sum_{\nu =0}^{j-1} a_{\nu} h_{\nu}(t) 
\end{equation*}
for some constants $a_{\nu}$, it is clear that any matrix of the form 
\begin{equation*}
(h_j(t_i)), \quad j, i = 0, 1, \ldots, m -1 , 
\end{equation*}
with all $t_i$ distinct can be transformed to a van der Monde matrix by elementary row operations, hence its determinant must be different from zero.

\smallskip
Lemma 3 shows that an arbitrary set of $m$ consecutive  rows 
$C_{2k} = (c_{2k, 0}, \ldots, c_{2k, m-1})$, $k = r, r+1, \ldots, r + m-1$,  
from the matrix $c_{2k, j}$ forms a linearly independent set of $m$-vectors. Therefore it is clear that an arbitrary row
$C_{2(r+m)}$ with $r \ge 0$ can be expressed as a linear combination of the $m$ preceding rows. This is made precise in the next lemma. 

\smallskip
\noindent
{\bf Lemma 4.}
{\sl 
The following identities hold: 
\begin{equation}       \label{lem2}
\sum_{k=0}^{m} (-1)^k \binom{m}k  c_{2r+2k ,j} = 0 \quad \tr{for  any $r \ge 0$ and $j = 0, 1, \ldots , m-1  $} . 
\end{equation}
}

\noindent
{\it Proof.}
Observe first of all that the identity
\begin{equation*}
\sum_{k=0}^{m} (-1)^k \binom{m}k h(k) = 0  
\end{equation*}
must be valid whenever $h(t)$ is a polynomial of degree at most $m-1$. This is obvious from the fact that the operator $h \mapsto \sum_{k=0}^{m} (-1)^k \binom{m}k h(k)$ is the composition of $m$ first order difference operators $h \mapsto h(t+1) - h(t)$.   As we saw above, the function  $k \mapsto c_{2k, j}$ is a polynomial of degree $j$, hence 
\begin{equation}   \label{ckj}
\sum_{k=0}^{m} (-1)^k \binom{m}k c_{2r+2k, j} = 0 \quad \tr{for $j \le m-1$} . 
\end{equation}
This completes the proof of Lemma 4.

\smallskip
\noindent
{\bf Corollary 2.}
{\sl 
Assume that the polynomials $p_{2k}$ are given by \eqref{id1}. Then the function $\rho(\w)^2$ must satisfy the identities 
\begin{align}   \label{pol-id-1}
&  \sum_{k=0}^{m} (-1)^k \binom{m}k \rho^{2(m-k)}  p_{2r + 2k} = 0 , \quad r = 0, 1, 2, \ldots , \quad \tr{or}  \\ 
 & p_{2r + 2m} = \sum_{k=0}^{m-1} r_k \, p_{2r + 2k} , 
 \label{rk-0}
\end{align}
where the coefficients $r_k$ are defined for $0 \le k \le m-1$ by 
\begin{equation}  
r_k  = - \binom mk (-\rho^2)^{m-k} = (-1)^{m-1} (-1)^k \binom mk \rho^{2(m-k)}  . 
\label{rk-2}
\end{equation}
}

\noindent
{\it Proof.}
Multiplying the respective equations in \eqref{id1} by suitable powers of $\rho^2$ and using Lemma 4 proves the assertion. 

For instance, if $m = 3$ the first few of the identities \eqref{rk-0} read 
\begin{align*}
p_6   & =  \rho^6 p_0  - 3 \rho^4 p_2 + 3 \rho^2 p_4           \\ 
p_8      & =   \rho^6 p_2  - 3 \rho^4 p_4 + 3 \rho^2 p_6           \\ 
p_{10}    & =  \rho^6 p_4  - 3 \rho^4 p_6 + 3 \rho^2 p_8           .   
 \end{align*}
The coefficients $r_k$ satisfy the identity 
\begin{equation}    \label{rk-1}
(t - \rho^2)^m  = t^m -  \sum_{k=0}^{m-1} r_k t^k . 
\end{equation}
And if we introduce the translation operator $T$,   defined by $T p_{2k} = p_{2k + 2}$,  on the infinite sequence of polynomials $(p_0, p_2, p_4, \ldots)$,  then \eqref{pol-id-1} can be written 
\begin{equation}      \label{T}
 (T - \rho^2)^{m} p_{2r} = 0 , \quad r = 0, 1, 2, \ldots . 
\end{equation}

\smallskip
A natural start towards a proof of  Theorem 2   would be to try to solve $\rho^2$ from some set of $m$ of the equations \eqref{pol-id-1} considering the equations as linear expressions in the $m$ unknowns $\rho^2, \ldots , \rho^{2 m}$.  If the matrix 
\begin{equation}     \label{A0}
A_0 = 
\begin{pmatrix}
 p_{0}       &     p_{2}        &  \ldots      &    p_{2m-2} \\ 
 p_{2}       &     p_{4}        &  \ldots      &    p_{2m}     \\   
\ldots        &        \ldots    &  \ldots      &      \ldots   \\ 
 p_{2m-2} &    p_{2m}      &  \ldots      &   p_{4m-4}  
\end{pmatrix}  
\end{equation} 
is non-singular, we can solve $\rho^2$ from the linear system \eqref{pol-id-1}  with $r = 0, 1, \ldots, m-1$ and obtain 
$\rho^2$ as a rational function 
\begin{equation*}
\rho^2 = F/G ,
\end{equation*}
where $F$ and $G$ are polynomials and $G = \det A_0$. As we shall see below (Lemma 5) it is easy to strengthen this argument by considering $m$ such systems together and thereby prove that $\rho^2$ is a polynomial. Therefore our main task in the rest of the proof of Theorem 2 will be to prove that the matrix $A_0$ is non-singular.

\smallskip
\noindent
{\bf Proposition 1.}
{\sl 
Let the polynomials $p_{2k}$ be defined as in Theorem 2 and assume that the function $q_{m-1}(\w)$ is not identically zero. Then the matrix \eqref{A0} is non-singular. 
}

The proof will be given at the end of this section.

Using Proposition 1 we can now easily finish the proof of Theorem 2. 
Denote by $\C(\w)$ the field of rational functions in $\w = (\w_1, \ldots, \w_n)$, and denote by $\C(\w)^m$ the $m$-dimensional vector space of $m$-tuples of elements from $\C(\w)$. 
Introduce the column $m$-vectors in $\C(\w)^m$ 
\begin{align*}
P_0  & = (p_0, p_2,  \ldots , p_{2m-2})^t, \quad P_2  = (p_2, p_4,  \ldots , p_{2m})^t, \\ 
       P_4  & = (p_4, p_6,  \ldots , p_{2m+2})^t, \ldots . 
\end{align*}
The recurrence relations \eqref{rk-0} then show that the translation operator $P_{2k} \mapsto P_{2k+2}$ is given by the matrix 
\begin{equation}     \label{S}
S = 
\begin{pmatrix}
0  & 1   &  0  & \ldots &  0   \\
 0  & 0   &  1  & \ldots &  0   \\
 \ldots  & \ldots   &  \ldots  & \ldots &  \ldots   \\
 0  & 0   &  \ldots  & 0 &  1   \\
 r_0  &  r_1   &  \ldots  & r_{m-2} &   r_{m-1}   
\end{pmatrix}  , 
\end{equation}
so that 
\begin{equation*}
S P_{2k} = P_{2k+2} \qquad   \tr{and} \quad  S^k P_0 = P_{2k} \quad \tr{for all $k$} . 
\end{equation*}
For instance, if $m = 3$ and $m = 4$, then 
\begin{equation*}
S = 
\begin{pmatrix}
0  & 1   &  0    \\
 0  & 0   &  1     \\
\rho^6  & - 3 \rho^4   &     3 \rho^2    
\end{pmatrix}
\qquad \tr{and} \quad
S = 
\begin{pmatrix}
0  & 1   &  0  &  0  \\
 0  & 0   &  1  &  0   \\  
  0  & 0   &  0  &  1   \\
- \rho^8  & 4 \rho^6   &  - 6 \rho^4   &  4 \rho^2 
\end{pmatrix}  , 
\end{equation*}
respectively. 
The characteristic equation of $S$ is 
\begin{align*}
\det(S - \lb I) & =  (-1)^{m-1}\big(r_0 - r_1 (-\lb) +  \ldots + (-1)^{m-1}(r_{m-1} - \lb) (-\lb)^{m-1} \big) \\ 
   & =  (-1)^{m-1}\big( r_0 + r_1 \lb + r_2 \lb^2 + \ldots + r_{m-1}\lb^{m-1} - \lb^m   \big)  . 
\end{align*} 
Defining $r_m = -1$ and taking \eqref{rk-1} into account we obtain 
\begin{equation*}
\det(S - \lb I) = \sum_{j=0}^m r_j \lb^j =  (\lb - \rho^2)^m  , 
\end{equation*}
so $S$ has the $m$-fold eigenvalue $\rho^2$, and the  
determinant of $S$ is 
\begin{equation*}
\det S =  (-1)^{m-1} r_0 = \rho^{2m} . 
\end{equation*}

\smallskip
\noindent
{\bf Lemma 5.}
{\sl 
Assume that the vectors $P_0, P_2 = S P_0, \ldots , P_{2m-2} = S^{m-1} P_0$ span $\C(\w)^m$, i.e., that the matrix \eqref{A0} is non-singular. Then $\rho(\w)^2$ is a polynomial. 
}

\noindent
{\it Proof.}  
We have already seen that $\rho(\w)^2$ must be a rational function, if the matrix \eqref{A0} is non-singular. 
The equations 
$P_{2k+2j} = S^k P_{2j}$ for $j = k, k+1, \ldots , k + m-1$ can be combined to the matrix equation   
\begin{equation}    \label{M}
A_k = S^k A_0 , 
\end{equation}
where $A_k$ is the matrix 
\begin{equation*}
A_k = 
\begin{pmatrix}
 p_{2k}       &     p_{2k+2}        &  \ldots      &    p_{2k+2m-2} \\ 
 p_{2k+2}       &     p_{2k+4}        &  \ldots      &    p_{2k+2m}     \\   
\ldots        &        \ldots    &  \ldots      &      \ldots   \\ 
 p_{2k+2m-2} &    p_{2k+2m}      &  \ldots      &   p_{2k+4m-4}  
\end{pmatrix}  . 
\end{equation*}
Denote the determinant of $A_0$, which is a polynomial, by $d(\w)$. Since $\det S = \rho^{2m}$, equation \eqref{M} implies that 
\begin{equation*}
\rho(\w)^{2mk} d(\w) 
\end{equation*}
is a polynomial for every $k$. Since $\rho^2$ is a rational function, this proves that $\rho^2$ must in fact be a polynomial as claimed.

\smallskip

We now turn to the proof of  Proposition 1. 
To motivate the next lemma we make the following observations. 
Let $B = (b_{k,j})$ be the matrix of the system \eqref{id1}, 
\begin{equation}    \label{bkj}
b_{k,j} = c_{2k,j} \rho^{2k-j} , \quad j = 0, 1, \ldots , m-1, \  k = 0, 1, \ldots ,  
\end{equation}
with $m$ columns and infinitely many rows, and let $B_0$ be the uppermost 
$m \times m$ minor of the matrix $B$ obtained by restricting $k$ to $0 \le k \le m-1$.  
Introducing the column $m$-vector $Q = (q_0, q_1, q_2, \ldots , q_{m-1})^t$ we then have 
$B_0 Q = P_0$. We want to prove that the vectors 
\begin{equation}   \label{span}
P_0, S P_0, \ldots , S^{m-1} P_0  \quad \tr{span} \quad \C(\w)^m .
\end{equation}
 Since $B_0 Q = P_0$, those vectors can be written 
$B_0 Q, S B_0 Q, \ldots , S^{m-1} B_0 Q$. And since $B_0$ is non-singular,    \eqref{span} is equivalent to 
\begin{equation*}  
(B_0^{-1} S B_0 )^{k} Q,   \quad k = 0, \ldots , m-1, \quad \tr{span} \quad \C(\w)^m .
\end{equation*}
Therefore we now study the matrix $B_0^{-1} S B_0$.

\smallskip
\noindent
{\bf Lemma 6.}
{\sl 
The matrix $B_0^{-1} S B_0$ is an upper triangular matrix of the form 
\begin{equation*}
B_0^{-1} S B_0 = \rho^2 I  +  N ,
\end{equation*}
where $N = (n_{k,j})$ is a nilpotent, upper triangular matrix whose elements next to the diagonal are given by 
\begin{equation}    \label{n}
n_{k, k+1} = 2 \rho \, k  , \quad 1 \le k \le m-1 . 
\end{equation}
}

For instance, if $m = 5$, 
\begin{equation*}
B_0^{-1} S B_0 = 
\begin{pmatrix}
\rho^2  &  2 \rho  &  2 &  0  &  0   \\  
0  &  \rho^2  &  4 \rho  &  6  &  0    \\  
0  &  0  &  \rho^2  &  6 \rho  &    12  \\  
0  &  0  &  0  &  \rho^2  &  8 \rho     \\  
0  &  0  &  0  &  0  &  \rho^2     
\end{pmatrix}  . 
\end{equation*}
The exact expression for the matrix $N$ is inessential for us, apart from the fact that all the entries  \eqref{n} are different from zero. 

\noindent
{\it Proof of Lemma 6.}
Denote by $u^0, \ldots, u^{m-1}$ the column vectors of $B_0$. 
The assertion of the lemma is that the matrix of $S$ with respect to the basis $u^0, \ldots , u^{m-1}$ is $\rho^2 I + N$. In fact we shall prove that 
\begin{align}
S u^0 & = \rho^2 u^0    \label{v1} \\  
S u^1 & = \rho^2 u^1  + 2 \, \rho \, u^0    \label{v2}\\  
S u^j & = \rho^2 u^j + 2 j \rho \, u^{j-1} + j (j-1)  \,  u^{j-2} ,    \label{v3}
     \quad 2 \le j \le m - 1 . 
\end{align}
The components of $u^j = (u^j_0, \ldots , u^j_{m-1})$ are 
\begin{equation*}
u^j_k = b_{k, j} = c_{2k,j}\, \rho^{2k-j} ,  \quad 0 \le j, \, k \le m-1 . 
\end{equation*}
Denote by $B_1$ the second uppermost $m \times m$ minor of the matrix $B$, which is obtained by 
restricting $k$ in \eqref{bkj} to $1 \le k \le m$. 
The argument in the proof of Corollary 2  shows that $S B_0 = B_1$, in other words
\begin{equation*}
S u^j = (u^j_1, \ldots ,  u^j_m) , \quad 0 \le j \le m-1 ,   
\end{equation*} 
if we define $u^j_m$ as $b_{m, j} = c_{2 m,j}\, \rho^{2 m - j}$. 
Denote by $D$ the formal derivative with respect to $\rho$. 
Note that $u^0_k = \rho^{2 k}$,  $u^1_k = D \rho^{2k}$ for $0 \le k \le m$, and that more  
  generally $u^j_k = D^{j} \rho^{2k}$ for $0 \le j \le m-1$, $0 \le k \le m$. 
The identity \eqref{v1} is obvious. To prove \eqref{v2} we just note that the $k$:th component of  $S u^1$ satisfies
\begin{equation*}
(S u^1)_k = D (\rho^2 \rho^{2k}) = \rho^2 D \rho^{2k} + 2  \rho \cdot \rho^{2k} = \rho^2 u^1_k + 2 \rho \, u^0_k , \quad  0 \le k \le m-1 . 
\end{equation*}
For \eqref{v3} we use Leibnitz' formula to get 
\begin{align*}
(S u^j)_k  & = D^{j} (\rho^2 \rho^{2k}) 
   = \rho^2 D^{j} \rho^{2k} 
                    + j 2 \rho  \, D^{j-1} \rho^{2k}  
                     + \binom j2 2  \, D^{j-2} \rho^{2k}   \\ 
  &  = \rho^2 u^j_k  +  2 j \rho \, u^{j-1}_k + j (j-1)  u^{j-2}_k  , \quad  2 \le j \le m-1 , \, 0 \le k \le m-1 , 
\end{align*}
which completes the proof.

\noindent
{\bf Lemma 7.}
{\sl
Let $A$ be an $m \times m$ matrix with entries in a  field $K$. Assume that $A$ has  one  $m$-fold eigenvalue $\lb$. Let $z\in K^m$.  Then the vectors 
$A^k z$, $k = 0, 1, \ldots , m-1$, span $K^m$ if and only if $(A - \lb I)^{m-1} z \ne 0$. 
}

\noindent
{\it Proof.}
The condition is obviously necessary, because if $(A - \lb I)^{m-1} z = 0$, then the equation 
\begin{equation*}
0 = (A - \lb I)^{m-1} z = \sum_{k=0}^{m-1} \binom {m-1}k A^k z (-\lb)^{m-1-k}  
\end{equation*}
shows that the $m$ vectors $A^k z$, $k = 0, 1, \ldots , m-1$,  are linearly dependent. 
The sufficiency follows from the fact that Jordan normal form (over the algebraic closure of $K$) must consist of just one Jordan block, but can also be seen more directly as follows. Introduce the subspace of $K^m$ defined by 
\begin{equation*}
E_k = \{x \in K^m;\, (A - \lb I)^{k} x = 0 \} \quad \tr{for} \quad  k = 1, \ldots , m .  
\end{equation*}
It is clear that $E_k \ss  E_{k+1}$, that $E_m = K^m$, and it is easily seen that $E_k = E_{k+1}$ implies $E_k = E_j$ for all $j \ge k$. 
The fact that $(A - \lb I)^{m-1} z \ne 0$ now shows that $0 \ne (A - \lb I)^k z \in E_{m-k} \sm E_{m-k-1}$ for $0 \le k \le m-1$, so the vectors 
$(A - \lb I)^k z $ must span $K^m$, which proves the lemma.

\noindent
{\it Proof of Proposition 1.}
By Lemma 7 it is enough to prove that $(B_0^{-1} S B_0 - \rho^2 I)^{m-1} Q \ne 0$.   
Using Lemma 6  we can   write 
\begin{equation*}
(B_0^{-1} S B_0 - \rho^2 I)^{m-1} Q  =   N^{m-1} Q   . 
\end{equation*} 
The matrix $N^{m-1}$ has only one non-vanishing entry in the upper right corner. The value of this entry is equal to  the product of of all the entries on the diagonal next to the main diagonal described in Lemma 6. And this product is equal to 
$c = (2 \rho)^{m-1} (m-1)!$, hence different from zero. Recalling that $Q = (q_0, q_1 , \ldots , q_{m-1})^t$ we can now conclude that 
\begin{equation*}
N^{m-1} Q = c \, q_{m-1}(\w) (1, 0, \ldots , 0)^t .  
\end{equation*}
By assumption the continuous function $q_{m-1}(\w)$ is different from zero on some open subset of $S^{n-1}$.  Thus we have proved  that the determinant of the matrix \eqref{A0} must be different from zero on some open set. But the determinant is a polynomial function, hence equal to a non-zero polynomial.  This completes the proof of Proposition 1.

\smallskip
\noindent
{\it End of proof of Theorem 1.}
Let $m$ be the largest number for which the coefficient $q_{m-1}$ in \eqref{g} is not identically zero. Then Proposition 1 shows that the matrix \eqref{A0} must be non-singular. And then Lemma 5 shows that $\rho(\w)^2$ must be a positive quadratic polynomial. This completes the proof of Theorem 2 and hence of Theorem 1.

\smallskip
\noindent
{\bf 5.  A semi-local result.}
The arguments given here prove in fact a semi-local version of  Theorem 1, where only an arbitrary open set of $\w$ comes into play, but all $p \in \R$. A set $W$ of hyperplanes $L \in \P$ is called translation invariant, if $L \in W$ implies that every translate $x + L$ is contained in $W$. 

\smallskip
\noindent
{\bf Theorem 3.}
{\sl Let $D$ be as in Theorem 1, let $x^0 \in \p D$, and let $\w^0$ be one of the unit normals of a supporting  plane $L_0$ to $\bar D$ at $x^0$.   If there exists a distribution $f$ with support in $\bar D$ and a translation invariant  open neighborhood $W$ of $L_0$, such that the restriction of  the distribution $Rf$ to $W$ is supported on the set of supporting planes to $D$ in $W$, then $\p D$ must be equal to the restriction of an ellipsoid in some neighborhood of $\pm x^0$. 
}

This theorem implies Theorem 1, because if the assumptions of Theorem 1 are fulfilled, then the assumptions of Theorem 3 are valid for every $x^0 \in \p D$, hence $\p D$ must be equal to an ellipsoid in some neighborhood of every point, hence be globally equal to an ellipsoid. 

\smallskip
\noindent
{\it Proof of Theorem 3.}
Let $f$  be as in the theorem and set $R f = g$. By the range conditions the functions 
$S^{n-1} \ni \w \mapsto \int_{\R} g(\w, p) p^k dp$ must belong to $\mathcal P_k$ for all  $k$.
Set $p^0 = x^0 \cdot \w^0$. The 
  assumptions imply that the restriction of $g$ to some neighborhood of $(\pm\w^0, \pm p^0)$ has the form \eqref{g}. This shows that the assumption \eqref{id1} of Theorem 2 must be valid for $\w$ in some neighborhood $E$ of $\pm\w^0$.  Note that the functions $q_j$ are defined only in some neighborhood of $\pm\w^0$, whereas $\rho(\w)$, the supporting function of $D$, is initially defined in all of $S^{n-1}$. But in this proof we are only concerned with the restriction of $\rho$ to $E$. Taking restriction to $E$ wherever relevant  we see that the proofs of Corollary 2,  Proposition 1, Lemma 5, and Lemma 6,  work without change, so the restriction of $\rho^2$ to  $E$ must be a positive quadratic polynomial,  and this implies the assertion.

    \end{document}